\documentclass[11pt,a4paper]{article}
\usepackage{fullpage}
\usepackage[cp1251]{inputenc}
\usepackage[russian, english]{babel}

\usepackage{amsthm}
\usepackage{amsmath}
\usepackage{amssymb}
\usepackage{amsfonts}
\usepackage{url,graphics}
\urldef\svetaurl\url{svetlana.obraztsova@gmail.com}
\newtheorem{theorem}{Theorem}[section]
\newtheorem{lemma}{Lemma}[section]
\newtheorem{arrow}{Corollary}[section]
\theoremstyle{definition}
\newtheorem{definition}{Definition}[section]
\newcommand{\bp}{\begin{proof}}
\newcommand{\ep}{\end{proof}}
\newcommand{\bl}{\begin{lemma}}
\newcommand{\el}{\end{lemma}}
\newcommand{\bt}{\begin{theorem}}
\newcommand{\et}{\end{theorem}}
\newcommand{\bd}{\begin{definition}}
\newcommand{\ed}{\end{definition}}
\newcommand{\ba}{\begin{arrow}}
\newcommand{\ea}{\end{arrow}}
\author{Svetlana Obraztsova\thanks{Steklov Institute of Mathematics at St.Petersburg, Russia and School of Physical and Mathematical Sciences, Nanyang Technological University, Singapore. Email: \svetaurl.}
\thanks{Supported by A*STAR SINGA Scholarship, partially supported by a grant for leading scientific schools support from the President of RF NSh-5282.2010.1, RFFI 09-01-12137-ofi-m}}
\title{On the chromatic uniqueness of $K_4$-homeomorphs with girth 7}
\begin{document}
\maketitle

\begin{abstract}
This paper settles the question left open in \cite{Peng3} (Discr. Math. 308~(2008), pp. 6132–-6140), completing the study of  $K_4$-homeomorphs i.e., cliques on $4$ vertices with edges replaced by paths, of girth 7. \end{abstract}

\section{Introduction}
This work contributes to a study of chromatic equivalence between
$K_4$ homeomorphs started in \cite{CZ} and continued
in series of articles. 

A thorough survey of the known results is presented in \cite{Koh1,Koh2}. The study of chromaticity of $K_4$-homeomorphs with at least three paths of same length has been completed in \cite{Ren} and  chromaticity of $K_4$-homeomorphs with at least two paths of length 1 also was presented in \cite{Xu} and \cite{Peng}. However, the chromatic uniqueness of some classes of $K_4$ homeomorphs remains open. In \cite{Peng2} the question of chromatic uniqueness of $K_4$ homeomorphs with girth at most 6 has been settled.

This work completes the study of chromaticity of $K_4$ homeomorphs with girth 7, which was initiated in \cite{Peng3}. The methods presented in this work can be applied to study the chromatic equivalence of graphs homeomorphic to $K_4$
with other girths, as well.

\section{Definitions}
We consider graphs which are non-directed, have no loops and multiple edges. We start with some definitions and notation. By $V(G)$ (respectively, by $E(G)$) we denote the set of vertices of graph $G$ (respectively, the set of edges) of graph $G=(V,E)$. Denote $n(G):=|V(G)|$ and $m(G):=|E(G)|$. By $K_4 (\alpha, \beta, \gamma, \delta, \epsilon, \eta)$ we denote the $K_4$-homeomorph, i.e. the graph K4 in which edges are replaced 
by paths with lengths $\alpha, \beta, \gamma, \delta, \epsilon, \eta$ (which will be called {\em parameters}), respectively (see fig.1). Let $C(G,k)$ denote the chromatic polynomial of $G$.
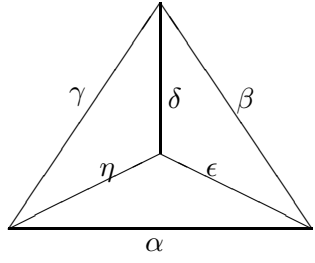
\begin{figure}[!h] \label{f:pict}
\begin{center}
\begin{picture}(50,40)
\put(5,5){\line(1,0){40}}
\put(5,5){\line(2,1){20}}
\put(45,5){\line(-2,1){20}}
\put(25,15){\line(0,1){20}}
\put(5,5){\line(2,3){20}}
\put(45,5){\line(-2,3){20}}
\put(23,2){$\alpha$}
\put(35,21){$\beta$}
\put(13,22){$\gamma$}
\put(26,21){$\delta$}
\put(17,12){$\eta$}
\put(31,12){$\epsilon$}
\end{picture}
\end{center}
\caption{$K_4( \alpha, \beta, \gamma, \delta, \epsilon, \eta)$}
\end{figure}

\bd
If $C(G,k)=C(J,k)$, then $G$ and $J$ are said to be chromatically equivalent and we denote this $G \sim J$.
\ed

\bd
A graph is chromatically unique, if $C(G,k)=C(J,k)$ implies that $G \cong J$. 
\ed

\section{Preliminaries}

In this section we survey some known results, which will be used in the sequel.

The following lemmas enable us to reduce the number of classes of chromatically  non-unique $K_4$-homeomorphs to be considered.
\bl[Koh and Teo \cite{Koh1}]\label{L:degree}
If $G \sim J$, then $|V(G)|=|V(J)|$ and $|E(G)|=|E(J)|$.
\el

\bl
[Chao and Zhao \cite{CZ}]Let $G$ and $J$ be undirected graphs without loops and multiple edges.
If $G \sim J$ and $G$ is homeomorphic to $K_4$, then $J$ is homeomorphic to $K_4$.
\el

\bl [Whitehead and Zhao \cite{WhZhao}] \label{L:length}
Let $K_4( \alpha, \beta, \gamma, \delta, \epsilon, \eta)\sim K_4( \alpha_1, \beta_1, \gamma_1, \delta_1, \epsilon_1, \eta_1)$. Then

\begin{enumerate}
	\item $min\{\alpha, \beta, \gamma, \delta, \epsilon, \eta\}=min\{\alpha_1, \beta_1, \gamma_1, \delta_1, \epsilon_1, \eta_1\}=:m^*$,
	\item the number of times that $m^*$ occurs in $(\alpha, \beta, \gamma, \delta, \epsilon, \eta)$ is equal to the number of times that $m^*$ occurs in  $(\alpha_1, \beta_1, \gamma_1, \delta_1, \epsilon_1, \eta_1)$.
\end{enumerate}
\el

\bl[Xu \cite{Xu}] \label{L:girth}
Assume that $K_4( \alpha, \beta, \gamma, \delta, \epsilon, \eta)\sim K_4( \alpha_1, \beta_1, \gamma_1, \delta_1, \epsilon_1, \eta_1)$. Then $K_4( \alpha, \beta, \gamma, \delta, \epsilon, \eta)$ and $K_4( \alpha_1, \beta_1, \gamma_1, \delta_1, \epsilon_1, \eta_1)$ have the same girth and the same number of cycles with the length equal to their girth.
\el

Next two lemmas give us the main method for studying chromatic uniqueness.

\bl[Li \cite{Li1} and Whitehead and Zhao \cite{WhZhao}] \label{L:polynomial}
The chromatic polynomial of $G=K_4( \alpha, \beta, \gamma, \delta, \epsilon, \eta)$ is
$$P(G,k)=(\frac{1}{k^2})(-1)^{m(G)}x[x^{m(G)-1}+Q(G,x)-(x+1)(x+2)],$$
where $x=1-k$ and $Q(G,x)=x^{\eta+\delta+\epsilon}+x^{\delta+\gamma+\beta}+x^{\alpha+\eta+\gamma}+x^{\alpha+\epsilon+\beta}+x^{\alpha+\delta}+x^{\eta+\beta}+x^{\gamma+\epsilon}-(x+1)(x^{\alpha}+x^{\beta}+x^{\gamma}+x^{\delta}+x^{\epsilon}+x^{\eta})$.

$Q(G,x)$ is called the essential polynomial of $G$.
\el

\bl[Li \cite{Li2}] 
Two $K_4$ homeomorphs with the same order are chromatically equivalent if and only if they have the same essential polynomial.
\el

The following lemma provides a convenient sufficient isomorphism condition of two $K_4$-homeomorphs.

\bl[Li \cite{Li1}] \label{l:iso}
Let $K_4( \alpha, \beta, \gamma, \delta, \epsilon, \eta)\sim K_4( \alpha_1, \beta_1, \gamma_1, \delta_1, \epsilon_1, \eta_1)$ and $\{\alpha, \beta, \gamma, \delta, \epsilon, \eta\}=\{ \alpha_1, \beta_1, \gamma_1, \delta_1, \epsilon_1, \eta_1 \}$ as multisets, then $K_4( \alpha, \beta, \gamma, \delta, \epsilon, \eta)\cong K_4( \alpha_1, \beta_1, \gamma_1, \delta_1, \epsilon_1, \eta_1)$.
\el

Next lemmas describe known results for chromatic uniqueness of graphs. These graphs can be subcases of $K_4$-homeomorphs with girth 7.

\bl[Ren \cite{Ren}]\label{l:ren} The $K_4$-homeomorph $K_4(a,b,c,d,e,f)$, where exactly three of $a,b,c,d,e,f$ are the same, is not chromatically unique if and only if $(a,b,c,d,e,f) \in \{(s,s,s-2,2,s,1), (s,s,s,1,s-2,2s-2), (t,t,1,t+2,t,2t), (t,t,1,t-1,t,2t), (t,t,t,1,t+1,2t+1), (1,1,1,3,t,t+1), (1,1,t,t+2,1,2)\} \quad s \geq 3, t \geq 2$.
\el

\bl[Peng and Liu \cite{Peng}] \label{l:peng2} The $K_4$-homeomorph $K_4(\alpha, 1, 1, \delta, \epsilon, \eta)$ ($\min{\alpha,\delta,\epsilon,\eta}>1$) is not chromatically unique if and only if $(\alpha, 1, 1, \delta, \epsilon, \eta) \in \{(a,1,1,a+b+1,b,b+1), (a,1,1,b,b+2,a+b), (a+1,1,1,a+3,2,a), (a+2,1,1,a,2,a +2), (3,1,1,2,b,b +1), (a+1,1,1,a,3,a+2), (a+1,1,1,b,3,a)\} \quad a \geq 2, b \geq 2$. Moreover,
$$K_4(a,1,1,a+b+1,b,b+1) \sim K_4(a,1,1,b,b+2,a+b),$$
$$K_4(a+1,1,1,a+3,2,a) \sim K_4(a+2,1,1,a,2,a+2),$$
$$K_4(1,b+2,b,1,2,2) \sim K_4(3,1,1,2,b,b+1),$$
$$K_4(1,a+1,a+3,1,2,a) \sim K_4(a+1,1,1,a,3,a+2),$$
$$K_4(1,a+2,b,1,2,a) \sim K_4(a+1,1,1,b,3,a).$$
\el
\bl[Xu \cite{Xu}] The $K_4$-homeomorph $K_4(1,\beta,\gamma,1,\epsilon,\eta)$ is not chromatically unique if and only if 
$(1,\beta,\gamma,1,\epsilon,\eta) \in \{(1,b+2,b,1,2,2), (1,a+1,a+3,1,2,a), (1, a + 2, b, 1, 2, a)\} \quad a \geq 2, b \geq 1$. Moreover,
$$K_4(1, b + 2, b, 1, 2, 2) \sim K_4(3, 1, 1, 2, b, b + 1),$$
$$K_4(1, a + 1, a + 3, 1, 2, a) \sim K_4(a + 1, 1, 1, a, 3, a + 2),$$
$$K_4(1, a + 2, b, 1, 2, a) \sim K_4(a + 1, 1, 1, b, 3, a).$$
\el

\bl [Peng \cite{Peng3}]\label{l:peng} The $K_4$-homeomorph $K_4(1,3,3,\delta,\epsilon,\eta)$, which has exactly one path of length 1 and has girth 7, is not chromatically unique if and only if $(1,3,3,\delta,\epsilon,\eta) \in\{(1,3,3,a-1,a,a+3),(1,3,3,a+1,a-1,a+2),(1,3,3,2,b,b+2),(1,3,3,2,4,7),(1,3,3,2,5,8),(1,3,3,5,2,5),(1,3,3,5,2,6),(1,2,4,3,7,3)\} \quad a>2,b \geq 2$. Moreover,
$$K_4(1,3,3,a-1,a,a+3) \sim K_4(1,3,3,a+1,a-1,a+2),$$
$$K_4(1,3,3,2,b,b+2) \sim K_4(1,2,4,b,b+1,3),$$
$$K_4(1,3,3,2,4,7)\sim K_4(1,2,4,4,3,6),$$
$$K_4(1,3,3,2,5,8)\sim K_4(1,2,4,6,3,6),$$
$$K_4(1,3,3,5,2,5)\sim K_4(1,2,4,3,3,6),$$
$$K_4(1,3,3,5,2,6)\sim K_4(1,2,4,3,7,3).$$

\el

The main method used in our work is different from methods that appear in the literature. First, we cancel out the common divisors when comparing the essential polynomials. Second, we replace $x$ (the unknown in the essential polynomials) by a root of the polynomial of $x^3+x+1$. This significantly reduces the number of subcases to consider.
\section{Main Result}

\bt \label{t:main}
Assume that the graph $G$ is a $K_4$-homeomorph with at most one path of length 1 and girth of $G$ is equal to 7.
Then $G$ is chromatically unique
unless it is one of the following:
$K_4(1,3,3,a-1,a,a+3)$,  $K_4(1,3,3,a+1,a-1,a+2)$,  $K_4(1,3,3,2,b,b+2)$, $K_4(1,2,4,b,b+1,3)$, $K_4(1,3,3,2,4,7)$, $K_4(1,2,4,4,3,6)$, $K_4(1,3,3,2,5,8)$, $K_4(1,2,4,6,3,6)$, $K_4(1,3,3,5,2,5)$, $K_4(1,2,4,3,3,6)$, $K_4(1,3,3,5,2,6)$, $K_4(1,2,4,3,7,3)$, $K_4(4,2,1,2,c+2,c)$, $K_4(3,2,2,c,1,c+3)$, $K_4(4,2,1,b,4,2)$, $K_4(2,2,3,b,5,1)$, $K_4(4,2,1,b,b+4,b+2)$, $K_4(4,2,1,b+1,b,b+5)$, $K_4(4,2,1,b+2,b,b+2)$, $K_4(4,2,1,b+1,b,b+3)$,  where $a>2, b \geq 2, c \geq 4, d \geq 5$.
\et 

A full list of chromatically equivalent graphs with girth 7 can be easily obtained by combining Lemmas \ref{l:ren}, \ref{l:peng2}, \ref{l:peng} and Theorem \ref{t:main}. 

\bp
We first list all classes of graphs, which could be $G$, and all classes of graphs, which could be chromatically equivalent to $G$.
From Lemmas \ref{L:girth} and \ref{L:length} it follows that these classes are the same. We now consider all possible cases for $G$: $K_4(3,3,1,c,b,a), K_4(1,a,2,2,2,b), K_4(4,2,1,c,b,a), K_4(2,2,3,c,b,a)$ (in last case $a$, $b$ or $c$ have to be equal to 1). Evidently, this list is exhaustive. The chromaticity of first and second classes is already settled in  Lemmas \ref{l:peng} and \ref{l:ren}, respectively. Thus, we need to consider only four cases.

\begin{eqnarray*} &&\text{case $(1)$: }K_4(4,2,1,\delta,\epsilon,\eta) \cong K_4(2,2,3,a,1,b), \\
&&\text{case $(2)$: }K_4(4,2,1,\delta,\epsilon,\eta) \cong K_4(3,2,2,b,a,1), \\
&&\text{case $(3)$: }K_4(4,2,1,\delta,\epsilon,\eta) \cong K_4(4,2,1,c,b,a), \\ 
&&\text{case $(4)$: }K_4(2,2,3,\delta,\epsilon,\eta) \cong K_4(2,2,3,a,b,c),
\end{eqnarray*}
where $\delta,\eta,\epsilon,a,b,c > 1.$ 
Denote a root of the polynomial $x^3+x+1$ by $t$.
Assume that $(1)$ happens. We need to find the parameters such that $K_4(4,2,1,\delta,\epsilon,\eta) \cong K_4(2,2,3,a,1,b)$. From \ref{L:polynomial} it follows that $K_4(4,2,1,\delta,\epsilon,\eta)\sim K_4(3,3,1,c,b,a)$ iff $Q(G,x)=Q(J,x)$. 

Let us write down this equality.
\begin{multline*}
x^{\eta+\delta+\epsilon}+x^{\delta+3}+x^{\eta+5}+x^{\epsilon+6}+x^{\delta+4}+x^{\eta+2}+x^{\epsilon+1}-(x+1)(x^{1}+x^{2}+x^{4}+x^{\delta}+x^{\epsilon}+x^{\eta})=\\ =x^{a+b+1}+x^{b+5}+x^{a+5}+x^{5}+x^{b+2}+x^{a+2}+x^{4}-(x+1)(x^{1}+x^{2}+x^{2}+x^{3}+x^{a}+x^{b}).
\end{multline*}

Since $\eta+\delta+\epsilon+7=a+b+8$ (see \ref{L:degree}), we see that $\eta+\delta+\epsilon=a+b+1$. Canceling same items in both sides and dividing by $x^2 -1$, we obtain

\begin{multline}\label{m:1}x^{\delta}(x^2+x+1)+(x^{\eta}-x^a-x^b)(x^3+x+1)+x^{\epsilon}(x^4+x^2+1)=2x^3+x^2.\end{multline}

The monomial $x^2$ occurs in RHS. Therefore, $x^2$ has to occur in RHS.
Hence, either
\begin{eqnarray*} &&\text{case $(1.1)$: }\delta=2,\\
&&\text{case $(1.2)$: }\epsilon=2, \text{or}\\
&&\text{case $(1.3)$: }\eta=2. \end{eqnarray*}

\begin{enumerate}
\item[(1.1)] If we replace $\delta$ by 2, we get the following identity after some simplification
\begin{multline}\label{m:11}(x^{\eta}-x^a-x^b)(x^3+x+1)+x^{\epsilon}(x^4+x^2+1)=x^3(1-x).\end{multline}
By substituting $t$ for $x$, we obtain $t^{\epsilon}(t^4+t^2+1)=t^3(1-t).$ It is easily seen that $t^4+t^2+1=1-t$. Thus, $t^{\epsilon}=t^3$ and $\epsilon=3$, because $t$ is neither a root of unity, nor zero. If we substitute 3 for $\epsilon$  in (\ref{m:11}), we obtain $(x^{\eta}-x^a-x^b)(x^3+x+1)=-x^4(x^3+x+1)$. This implies that w.l.o.g. $\eta=b$ and $a=4$. Clearly, by Lemma \ref{l:iso} we get a pair of isomorphic graphs $K_4(4,2,1,2,3,b) \cong K_4(2,2,3,4,1,b)$. Thus, case $(1.1)$ cannot happen.
\end{enumerate}

\item[(1.2)] By substituting 2 for $\epsilon$ in (\ref{m:1}), we get 
\begin{multline}\label{m:12}x^{\delta}(x^2+x+1)+(x^{\eta}-x^a-x^b)(x^3+x+1)=x^3(2-x-x^3).\end{multline}
The monomial $2x^3$ occurs in RHS, thus, $2x^3$ has to occur in LHS. Therefore, we get either $\delta=\eta=3$, $\delta=3,\eta=2$, or $\delta=2,\eta=3$. It can be easily checked that if $\delta=2$ or $\eta=2$, we have w.l.o.g. $a=2$, because the monomial $x^2$ has to be canceled on LHS. This case, $x^3$ cancels on RHS. That is $\delta=\eta=3$.
If we replace $\delta$ and $\eta$ by 3 and divide by $x^4$, we obtain
$$-(x^{a-4}+x^{b-4})(x^3+x+1)=-(x^2+x+2).$$ Evidently, product of roots of RHS equals 2 and product of roots of LHS equals either 1, or 0. Hence, case $(1.2)$ cannot happen.

\item[(1.3)] If we substitute 2 for $\eta$ in (\ref{m:1}), we obtain
\begin{multline}\label{m:13}x^{\delta}(x^2+x+1)-(x^a+x^b)(x^3+x+1)+x^{\epsilon}(x^4+x^2+1)=x^3(1-x^2).\end{multline}
The monomial $x^3$ has to occur in LHS. Thus, since cases $\delta=2$ and $\eta=2$ have been already considered above,we have either $\delta=3$ (case $(1.3.1)$), or $\epsilon=3$ (case $(1.3.2)$). 

\item[(1.3.1)] By replacing $\delta$ by 3 in (\ref{m:13}), we get
\begin{multline}\label{m:131}-(x^a+x^b)(x^3+x+1)+x^{\epsilon}(x^4+x^2+1)=-x^4(1+2x).\end{multline}
By substituting $t$ for $x$, we obtain $(t^{\epsilon}(t^4+t^2+1)=-t^4(1+2t).$ We see that $t^4+t^2+1=1-t$ and $1+2t=t-t^3=t(1-t^2).$ It follows that $t^{\epsilon}=t^8$. Hence, $\epsilon=8$. If we substitute 8 for $\epsilon$ in (\ref{m:131}) and divide by $-(x^3+x+1)$, we obtain $x^a+x^b=x^9-x^6+x^5+x^4.$ Thus, the case $(1.3.1)$ cannot happen.

\item[(1.3.2)] If we substitute 3 for $\epsilon$ in (\ref{m:13}), we obtain
\begin{multline*}x^{\delta}(x^2+x+1)-(x^a+x^b)(x^3+x+1)=-x^5(x^2+2).\end{multline*}
By substituting $t$ for $x$, we have $t^{\delta}(t^2+t+1)=-t^5(t^2+2).$ It is clear that $t^2+t+1=1-t$ and $-t(2+t^2)=1-t$. It follows that $t^{\delta+2}=t^4$ and $\delta=2$. This case has been already considered, therefore, the case $(1.3.2)$ cannot happen.

We looked through the case $(1)$ and this case is impossible.
The techniques applied to this case can be similarly applied to the remaining cases (namely, $(2)$, $(3)$ and $(4)$).
We omit these cases
(and their numerous subcases) here, due to their technical similarity; the reader can find 
these cases in Appendix A.

\ep

\section{Acknowledgments}
My sincere thanks to Dmitrii V Pasechnik for his great help with all aspects of the paper.

\clearpage
\normalsize
\section*{Appendix A}

We consider here cases $(2)$, $(3)$ and $(4)$ in Theorem~\ref{t:main}.

\begin{itemize}
\item[(2)]Consider the equality
\begin{multline*}
x^{\delta+3}+x^{\eta+5}+x^{\epsilon+6}+x^{\delta+4}+x^{\eta+2}+x^{\epsilon+1}-(x+1)(x^{4}+x^{\delta}+x^{\epsilon}+x^{\eta})=\\ =x^{b+4}+x^{6}+x^{a+5}+x^{b+3}+x^{a+2}+x^{3}-(x+1)(x^{2}+x^{3}+x^{a}+x^{b}).
\end{multline*}
If we divide by $x+1$, cancel $-x^4$ and divide by $x-1$, we obtain
\begin{multline}\label{m:2}
(x^{\delta}-x^b)(x^2+x+1)+(x^{\eta}-x^a)(x^3+x+1)+x^{\epsilon}(x^4+x^2+1)=x^2(x^2+x+1).
\end{multline}
The monomial $x^2$ occurs in RHS, thus, $x^2$ has to occur in LHS. Therefore,
we have one of the following
\begin{eqnarray*} 
&&\text{case $(2.1)$: }\delta=2, \\
&&\text{case $(2.2)$: }\epsilon=2,\\
&&\text{case $(2.3)$: }\eta=2.\end{eqnarray*}
\item[(2.1)] If we subsitute 2 for $\delta$ in (\ref{m:2}), we get 
\begin{multline}\label{m:21}
-x^b(x^2+x+1)+(x^{\eta}-x^a)(x^3+x+1)+x^{\epsilon}(x^4+x^2+1)=0.
\end{multline}
By substituting $t$ for $x$, we have $t^b(t^2+t+1)=t^{\epsilon}(t^4+t^2+1).$ It is easily shown that $t^2+t+1=t^2-t^3$ and $t^4+t^2+1=1-t$, whence, $t^{b+2}=t^{\epsilon}$. This implies that $\epsilon=b+2$. That is, we have the pair of chromatically equivalent graphs $K_4(4,2,1,2,b+2,b) \cong K_4(3,2,2,b,b+3,1),$ where $b>2$.
\item[(2.2)] If we replace $\epsilon$ by 2 in (\ref{m:2}), we obtain
\begin{multline}\label{m:22}
(x^{\delta}-x^b)(x^2+x+1)+(x^{\eta}-x^a)(x^3+x+1)=x^3(1-x^3).
\end{multline}
The monomial $x^3$ occurs in RHS, thus, there is term in LHS, which is equal to $x^3$. Therefore, either
\begin{eqnarray*} &&\text{case $(2.2.1)$: } \delta=3, \\
&&\text{case $(2.2.2)$: } \eta=2,  \text{or}\\
&&\text{case $(2.2.3)$: }\eta=3.\end{eqnarray*}
\item[(2.2.1)] By substituting 3 for $\delta$ in (\ref{m:22}), we get
\begin{multline*}
-x^b(x^2+x+1)+(x^{\eta}-x^a)(x^3+x+1)=-x^4(x^2+x+1).
\end{multline*}
If we replace $x$ by $t$, we get $-t^b(t^2+t+1)=-t^4(t^2+t+1)$. Hence, $t^b=t^4$ and $b=4$. Therefore, $(x^{\eta}-x^a)(x^3+x+1)=0$ and $\eta=a$. Obviously, by Lemma \ref{l:iso} we get a pair of isomorphic graphs $K_4(4,2,1,3,2,a) \cong K_4(3,2,2,4,a,1)$. Thus, case $(2.2.1)$ cannot happen.
\item[(2.2.2)] If we sustitute 2 for $\eta$ in (\ref{m:22}), we obtain
\begin{multline*}
(x^{\delta}-x^b)(x^2+x+1)-x^a(x^3+x+1)=-x^2(x^4+x^3+1).
\end{multline*} 
The monomial $-x^2$ occurs in RHS, thus, $-x^2$ has to occur in LHS and either $a=2$, or $b=2$. This case, 3 parameters of $K_4(3,2,2,b,a,1)$ are equal to 2. Using lemma \ref{l:ren} it can easily be proved that $K_4(3,2,2,b,a,1)$ is chromatically unique. That is, case $(2.2.2)$ cannot happen.
\item[(2.2.3)] By replacing $\eta$ by 3 in (\ref{m:22}), we have
\begin{multline}\label{m:223}
(x^{\delta}-x^b)(x^2+x+1)-x^a(x^3+x+1)=-x^4(1+2x^2).
\end{multline}
The monomial $-x^4$ occurs in RHS, hence, one of the terms in LHS has to be equal to $-x^4$. Therefore, either
\begin{eqnarray*} &&\text{case $(2.2.3.1)$: }b=3, \\
&&\text{case $(2.2.3.2)$: }b=4,\\
&&\text{case $(2.2.3.3)$: }a=3, \text{ or}\\
&&\text{case $(2.2.3.4)$: }a=4. \end{eqnarray*}
We can skip case $b=2$, because it considers as case $(2.2.2)$.
\item[(2.2.3.1)] If we substitute 4 for $b$ in (\ref{m:223}), we get
\begin{multline*}
x^{\delta}(x^2+x+1)-x^a(x^3+x+1)=x^3(1+x^2-2x^3).
\end{multline*} 
There is term $x^3$ in RHS, therefore, the monomial $x^3$ has to occur in LHS. Since the case of $\delta=2$ has been discussed in case $(2.1)$, we can suppose $\delta=3$. This case we have $x^a(x^3+x+1)=-x^4(1+2x^2).$ Evidently,this equality does not hold. Hence, case $(2.2.3.1)$ cannot happen.

\item[(2.2.3.2)] If we replace $b$ by 5 in (\ref{m:223}), we obtain
\begin{multline*}
x^{\delta}(x^2+x+1)-x^a(x^3+x+1)=x^5(1-x).
\end{multline*} 
If we substitute $t$ for $x$, we have $t^{\delta}(t^2+t+1)=t^5(1-t).$ Clearly, $t^2+t+1=t^2-t^3$, hence, $t^{\delta+2}=t^5$ and $\delta=3$. By substituting 3 for $\delta$, we get $-x^a(x^3+x+1)=-x^3(x^3+x+1)$ and $a=3$. Thus, we obtain chromatically equivalent graphs $K_4(4,2,1,3,5,3)$ and $K_4(3,2,2,5,3,1)$. Therefore, $G \cong J$ by using Lemma \ref{l:iso}. That is, the case $(2.2.3.2)$ cannot happen.

\item[(2.2.3.3)] By substituting 4 for $a$ in (\ref{m:223}), we get $(x^{\delta}-x^b)(x^2+x+1)=x^3(1-x^3)$. Whence, $x^{\delta}-x^b=x^3(1-x)$ and $\delta=3$, $b=4$. This case cannot happen similar to case $(2.2.3.2)$. 

\item[(2.2.3.4)] By replacing $a$ by 5 in (\ref{m:223}), we get $(x^{\delta}-x^b)(x^2+x+1)=x^5(x-1)^2$. Evidently, LHS cannot be divided by $(x-1)^2$. Thus, the case $(2.2.3.4)$ cannot happen.

\item[(2.3)] If substitute 2 for $\eta$ in (\ref{m:2}), we obtain
\begin{multline}\label{m:23}
(x^{\delta}-x^b)(x^2+x+1)-x^a(x^3+x+1)+x^{\epsilon}(x^4+x^2+1)=x^4(1-x).
\end{multline}
The monomial $x^4$ occurs in RHS, thus, there is term in LHS which are equal to $x^4$. Whence, we have one of the following
\begin{eqnarray*} &&\text{case $(2.3.1)$: }\delta=3, \\
&&\text{case $(2.3.2)$: } \delta=4,\\
&&\text{case $(2.3.3)$: }\epsilon=4. \end{eqnarray*}
We can skip cases $b=2$ and $\epsilon=2$, because we consider these cases as $(2.1)$ and $(2.2)$ respectively.
\item[(2.3.1)] If we substitute 3 for $\delta$ in (\ref{m:23}), we get
\begin{multline}\label{m:231}
-x^b(x^2+x+1)-x^a(x^3+x+1)+x^{\epsilon}(x^4+x^2+1)=-x^3(1+2x^2).
\end{multline}
If $a=2$ or $b=2$ than $K_4(3,2,2,b,a,1)$ has 3 parameters, which equal to 2. This case, using Lemma \ref{l:ren} it follows easily that $K_4(3,2,2,b,a,1)$ is chromatically unique. On the other hand, the monomial $-x^3$ occurs in RHS, therefore, one term of LHS has to be equal to $-x^3$. This implies that either
\begin{eqnarray*}
&&\text{case $(2.3.1.1)$: } b=3, \text{ or} \\
&&\text{case $(2.3.1.2)$: } a=3.
\end{eqnarray*}
\item[(2.3.1.1)] By replacing  $b$ by 3 in (\ref{m:231}), we obtain $-x^a(x^3+x+1)+x^{\epsilon}(x^4+x^2+1)=-x^4(1-x).$
By substituting $t$ for $x$, we have $t^{\epsilon}(t^4+t^2+1)=-t^4(1-t)$. Trivially, $t^4+t^2+1=1-t$, whence, $t^{\epsilon}=t^4$ and $\epsilon=4$. It follows that $-x^a(x^3+x+1)=-x^5(1+x+x^3)$ and $a=5$. In other words, we obtain a pair of chromatically equivalent graphs $K_4(4,2,1,3,4,2) \cong K_4(3,2,2,3,5,1).$

\item[(2.3.1.2)] If we substitute 3 for $a$ in (\ref{m:231}), we get $-x^b(x^2+x+1)+x^{\epsilon}(x^4+x^2+1)=x^4(1-x)^2.$ Obviously, LHS can be divided by $x^2+x+1$ and RHS cannot. Thus, this case cannot happen.

\item[(2.3.2)] If we replace $\delta$ by 4 in (\ref{m:23}), we obtain 
\begin{multline}\label{m:232}
-x^b(x^2+x+1)-x^a(x^3+x+1)+x^{\epsilon}(x^4+x^2+1)=-x^5(2+x).
\end{multline}
It is obvious that $x^{\epsilon+4}$ has to be canceled on LHS, hence, either 
\begin{eqnarray*}
&&\text{case $(2.3.2.1)$: }\epsilon+2=b, \text{ or} \\ 
&&\text{case $(2.3.2.2)$: }\epsilon+1=a 
\end{eqnarray*}
(otherwise $\deg LHS > \epsilon+4 \geq 6=\deg RHS$).

\item[(2.3.2.1)] By substituting $\epsilon+2$ for $b$ and $t$ for $x$ in (\ref{m:232}), we get $t^{\epsilon}(1-t^3)=-t^5(2+t)$. It is easy to see that $2+t=1-t^3$. Therefore, $t^{\epsilon}=-t^5$. Such $\epsilon$ does not exist, the case $(2.3.2.1)$ cannot happen.

\item[(2.3.2.2)] If we replace $a$ by $\epsilon+1$ in (\ref{m:232}), we obtain $-x^b(x^2+x+1)+x^{\epsilon}(1-x)=-x^5(2+x)$. The monomial $-2x^5$ occurs in RHS, thus, two terms in LHS have to be equal to $-x^5$. Thus, $\epsilon=4$. Therefore, $-x^b(x^2+x+1)=-x^4(1+x+x^2)$ and $b=4$. That is, we have a pair of chromatically equivalent graphs $K_4(4,2,1,4,4,2) \cong K_4(3,2,2,4,5,1).$

\item[(2.3.3)] By substituting 4 for $\epsilon$ and $t$ for $x$ in (\ref{m:232}), we have $(t^{\delta}-t^b)(t^2+t+1)=0$. Thus, $\delta=b$. This means that we have a pair of chromatically equivalent graphs $K_4(4,2,1,b,4,2) \cong K_4(3,2,2,b,5,1),$ where $b \geq 3$. Pairs of graphs, which we obtain in cases $(2.3.2.2)$ and $(2.3.2.2)$, are this pair with $b=3$ and $b=4$ respectively.

That is, we looked trough case $(2)$ and case $(2)$ happens iff either $K_4(4,2,1,2,c+2,c) \cong K_4(3,2,2,c,1,c+3)$, or $K_4(4,2,1,b,4,2)\cong K_4(2,2,3,b,5,1)$.

\item[(3)] If we divide by $x^2-1$ both sides of equality for essential polynomial, we get
\begin{multline}\label{m:3}x^{\delta}(x^2+x+1)+x^{\eta}(x^3+x+1)+x^{\epsilon}(x^4+x^2+1)=x^{c}(x^2+x+1)+x^{a}(x^3+x+1)+x^{b}(x^4+x^2+1).\end{multline}
Evidently, $\min\{\delta,\epsilon,\eta\}=\min\{a,b,c\}.$ Therefore, we have 6 possibilities, w.l.o.g. we have one of the following
\begin{eqnarray*}
&&\text{case $(3.1)$: } \delta=a=\min\{\delta,\epsilon,\eta\}=\min\{a,b,c\},\\
&&\text{case $(3.2)$: }\delta=b=\min\{\delta,\epsilon,\eta\}=\min\{a,b,c\}, \delta <c, \\
&&\text{case $(3.3)$: }\delta=c=\min\{\delta,\epsilon,\eta\}=\min\{a,b,c\}, \delta <a, \delta <c, \\
&&\text{case $(3.4)$: }\eta=a=\min\{\delta,\epsilon,\eta\}=\min\{a,b,c\}, \eta <c, \\
&&\text{case $(3.5)$: }\eta=b=\min\{\delta,\epsilon,\eta\}=\min\{a,b,c\}, \eta <c, \eta <a, \text{ or}\\
&&\text{case $(3.6)$: }\epsilon=a=\min\{\delta,\epsilon,\eta\}=\min\{a,b,c\}, \epsilon <a,c.\end{eqnarray*}

\item[(3.1)] If we substitute $c$ for $\delta$ in (\ref{m:3}), cancel the same items in both sides and replace $x$ by $t$, we have $t^{\epsilon}(t^4+t^2+1)=t^{b}(t^4+t^2+1).$ Therefore, $t^{\epsilon}=t^{b}$ and $\epsilon=b$. Hence, $\eta=a$. In other words, we obtain a pair of isomorphic graphs.
\item[(3.2)] If we replace $a$ by $\delta$ in (\ref{m:3}) and divide by $\delta$ both sides of equality, we get
\begin{multline}\label{m:32}x^{\eta-\delta}(x^3+x+1)-x^{c-\delta}(x^2+x+1)+(x^{\epsilon-\delta}-x^{b-\delta})(x^4+x^2+1)=x^2(x-1).\end{multline}
The monomial $-x^2$ occurs in RHS, thus, $-x^2$ has to occur in LHS. Since $c>\delta$ this implies that either
\begin{eqnarray*}
&&\text{case $(3.2.1)$: } c-\delta=1, \\
&&\text{case $(3.2.2)$: } c-\delta=2,  \\
&&\text{case $(3.2.3)$: } b-\delta=2, \text{ or}\\
&&\text{case $(3.2.4)$: } b-\delta=0.\end{eqnarray*}
\item[(3.2.1)] By substituting 1 for $c-\delta$ in (\ref{m:32}), we have $x^{\eta-\delta}(x^3+x+1)+(x^{\epsilon-\delta}-x^{b-\delta})(x^4+x^2+1)=x(1+2x^2).$ Clearly, one term of LHS has to be equal to $x$, hence, either $\eta-\delta=1$, or $\epsilon-\delta=1$. If $\eta-\delta=1$ than $(x^{\epsilon-\delta}-x^{b-\delta})(x^4+x^2+1)=-x^2(1-x)^2.$ We see that RHS can be divided by $(1-x)^2$ and RHS cannot. Thus, $\epsilon-\delta=1$ and $x^{\eta-\delta}(x^3+x+1)-x^{b-\delta}(x^4+x^2+1)=x^3(1-x^2).$ If we replace $x$ by $t$, we have $-t^{b-\delta}(t^4+t^2+1)=t^3(1-t^2)=-t^6(1-t)$. Since $t^4+t^2+1=1-t$ we get $-t^{b-\delta}=-t^6$ and $b-\delta=6$. Therefore, $x^{\eta-\delta}(x^3+x+1)=x^3-x^5+x^{10}+x^{8}+x^{6}.$ Evidently, this equality does not hold. That is, the case $(3.2.1)$ cannot happen.
\item[(3.2.2)] If we replace $c-\delta$ by 2 in (\ref{m:32}), we get
$x^{\eta-\delta}(x^3+x+1)+(x^{\epsilon-\delta}-x^{b-\delta})(x^4+x^2+1)=x^3(x+2).$
There is the monomial $2x^3$ in RHS, thus, two terms of LHS have to be equal to $x^3$. This yields that either $\eta-\delta=2$, or $\eta-\delta=3$ and either $\epsilon-\delta=1$, or $\epsilon-\delta=3$. If $\eta-\delta=2$ than $(x^{\epsilon-\delta}-x^{b-\delta})(x^4+x^2+1)=-x^2(1+x)(1-x)^2.$ Evidently, LHS cannot be divided by $(1-x)^2$, Thus, identity does not hold. If $\eta-\delta=2$ than $(x^{\epsilon-\delta}-x^{b-\delta})(x^4+x^2+1)=x^3(1-x^3).$ Evidently, RHS can be divided by $x^2+x+1$ and LHS cannot, thus, identity does not hold. This means that case $(3.2.2)$ cannot happen.
\item[(3.2.3)] By substituting 2 for $b-\delta$ in (\ref{m:32}), we have
$x^{\eta-\delta}(x^3+x+1)-x^{c-\delta}(x^2+x+1)+x^{\epsilon-\delta}(x^4+x^2+1)=x^3(x^3+x+1).$ 
If we substitute $t$ for $x$, we obtain $t^{\epsilon-\delta}(t^4+t^2+1)=t^{c-\delta}(t^2+t+1).$ Since $t^4+t^2+1=1-t$ and $t^2+t+1=t^2-t^3$ we have $t^{\epsilon-\delta}=t^{c-\delta+2}$ and $\epsilon-\delta=c-\delta+2$. If we replace $\epsilon -\delta$ by $c -\delta+2$, we get $x^{\eta-\delta}(x^3+x+1)+x^{c-\delta}(x^6+x^4-x-1)=x^3(x^3+x+1).$ Dividing identity by $x^3+x+1$, we have $x^{\eta-\delta}+x^{c-\delta}(x^3-1)=x^3.$ By our assumption $c>\delta$, thus, $\deg LHS \geq 4 > \deg RHS$, therefore, this identity does not hold. That is, the case $(3.2.3)$ cannot happen.

\item[(3.2.4)] If we replace $b-\delta$ by 0 in (\ref{m:32}), we get 
$x^{\eta-\delta}(x^3+x+1)-x^{c-\delta}(x^2+x+1)+x^{\epsilon-\delta}(x^4+x^2+1)=x^4+x^3+1.$
There is term $1$ in RHS, hence, one term in LHS has to be equal to 1. If $\epsilon-\delta=0$ than $x^{\eta-\delta}(x^3+x+1)-x^{c-\delta}(x^2+x+1)=x^3-x^2$ and either $\eta-\delta=1$, or $\eta-\delta=2$. These cases have been already considered as cases $(3.2.1)$ and $(3.2.2)$ respectively. Whence, $\eta-\delta=0$ and we have $-x^{c-\delta}(x^2+x+1)+x^{\epsilon-\delta}(x^4+x^2+1)=x(x^3-1)$. Dividing both sides of this equality by $x^2+x+1$, we get $x^{\epsilon-\delta}(x^2-x+1)-x^{c-\delta}=x(x-1)$. Clearly, $x^{\epsilon-\delta+1}$ cannot be canceled on LHS, thus, $\epsilon-\delta=0$ and $-x^{c-\delta}=-1$. Since $c>\delta$ this identity does not hold. 

\item[(3.3)] By substituting  $\delta$ for $b$ in (\ref{m:3}) and dividing the equality by $x^{\delta}$, we obtain
\begin{multline}\label{m:33}(x^{\eta-\delta}-x^{a-\delta})(x^3+x+1)+x^{\epsilon-\delta}(x^4+x^2+1)-x^{c-\delta}(x^2+x+1)=x(x^3-1).\end{multline}
The monomial $-x$ occurs in LHS, therefore, one term in RHS is equal to $-x$. Since $\delta <a$ and $\delta<c$ than either
\begin{eqnarray*}
&&\text{case (3.3.1): } a=\delta+1, \text{ or} \\
&&\text{case (3.3.2): } c=\delta+1.
\end{eqnarray*}
\item[(3.3.1)] If we substitute $\delta+1$ for $a$ in (\ref{m:33}), we get
\begin{multline}\label{m:331}x^{\eta-\delta}(x^3+x+1)+x^{\epsilon-\delta}(x^4+x^2+1)-x^{c-\delta}(x^2+x+1)=x^2(2x^2+1).\end{multline} 
There is term $x^2$ in RHS, thus, $x^2$ has to occur in LHS. Hence, either
\begin{eqnarray*} &&\text{case $(3.3.1.1)$: }\epsilon=\delta, \\
&&\text{case $(3.3.1.2)$: } \epsilon=\delta+2, \\
&&\text{case $(3.3.1.3)$: } \eta=\delta+2, \text{or}\\
&&\text{case $(3.3.1.3)$: } \eta=\delta+1. \end{eqnarray*} 

\item[(3.3.1.1)] If we replace $\epsilon$ by $\delta$ and $x$ by $t$, in (\ref{m:331}), we obtain 
$-t^{c-\delta}(t^2+t+1)=t^4-1$. Since $t^2+t+1=t^2-t^3$ we have $t^{c - \delta+2}=1+t+t^2+t^3=t^2$. By our assumption $c>\delta$, hence, this equality does not hold.

\item[(3.3.1.2)] If we substitute $\delta+2$ for $\epsilon$ in (\ref{m:331}), we obtain
$x^{\eta-\delta}(x^3+x+1)-x^{c-\delta}(x^2+x+1)=x^4(1-x^2).$
By replacing $x$ by $t$, we get 
$-t^{c-\delta}(t^2+t+1)=t^4(1+t)(1-t)$. Since $t^2+t+1=1-t$ and $1+t=-t^3$ we have $t^{c-\delta+2}=t^7$ and $c=\delta+5$. If we substitute $\delta+5$ for $c$, we get $x^{\eta-\delta}(x^3+x+1)=x^4(1+x+x^3)$ and $\eta=\delta+4$.
That is, we obtain a pair of chromatically equivalent graphs $K_4(4,2,1,\delta,\delta+2,\delta+4) \cong K_4(4,2,1,\delta+1,\delta,\delta+5),$ where $\delta>1$.

\item[(3.3.1.3)] By substituting $\delta+2$ for $\eta$ in (\ref{m:331}), we get
$x^{\epsilon-\delta}(x^4+x^2+1)-x^{c-\delta}(x^2+x+1)=-x^3(1-x)^2$. It is easily seen that LHS can be divided by $x^2+x+1$ and RHS cannot. In other words, the case $(3.3.1.3)$ cannot happen.

\item[(3.3.1.4)] If we substitute $\delta+1$ for $\eta$ in (\ref{m:331}) and divide both sides of equality by $x^2+x+1$, we obtain
$x^{\epsilon-\delta}(x^2-x+1)-x^{c-\delta}=x^2-x$.
There are 1 positive monomial in RHS, whence, either $x^{\epsilon-\delta+2}$, or $x^{\epsilon-\delta}$ has to be canceled on LHS. The degree of positive term in LHS is greater than the degree of negative term in LHS. It implies that $x^{\epsilon-\delta}$ has to be canceled and $\epsilon-\delta=c-\delta$. Equivalently, we have pair of isomorphic graphs $K_4(4,2,1,\delta,\delta,\delta+1) \cong K_4(4,2,1,\delta,\delta,\delta+1)$. Thus, cases $(3.3.1.4)$ and $(3.3.1)$ do not happen.

\item[(3.3.2)] If we replace $c$ by $\delta+1$ in (\ref{m:33}), we get
$(x^{\eta-\delta}-x^{a-\delta})(x^3+x+1)+x^{\epsilon-\delta}(x^4+x^2+1)=x^2(x^2+x+1).$
By substituting $t$ for $x$, we obtain $t^{\epsilon-\delta}(t^4+t^2+1)=t^2(t^2+t+1).$ Since $t^4+t^2+1=1-t$ and $t^2+t+1=t^2-t^3$ we have $t^{\epsilon-\delta}=t^4$ and $\epsilon=\delta+4$. If we substitute $\delta+4$ for $\epsilon$ and divide the equality by $x^3+x+1$, we get $x^{\eta-\delta}-x^{a-\delta}=x^2(1-x^3).$ Therefore, $\eta-\delta=2$ and $a-\delta=5$. That is, we have a pair of chromatically equivalent graphs $K_4(4,2,1,\delta,\delta+2,\delta+4) \cong K_4(4,2,1,\delta+1,\delta,\delta+5)$, but this pair was already found in case $(3.3.1.2)$.

\item[(3.4)] By substituting $a$ for $\eta$ in (\ref{m:3}) and dividing by $x^2+x+1$, we get $x^{\delta}+x^{\epsilon}(x^2-x+1)=x^{c}+x^{b}(x^2-x+1).$
Remember, $\delta+\epsilon=b+c$. Obviously, minimal degree in RHS is equal to either $b$, or $c$ and minimal degree in LHS is equal to either $\delta$, or $\epsilon$. Hence, either $\delta=b$ and $\epsilon=c$, or $\delta=c$ and $\epsilon=b$. Both cases we have a pair of isomorphic graphs by lemma (\ref{l:iso}).

\item[(3.5)] By substituting $\eta$ for $b$ in (\ref{m:3}) and dividing by $x^b$, we obtain
\begin{multline}\label{m:35}(x^{\delta-b}-x^{c-b})(x^2+x+1)+x^{\epsilon-b}(x^4+x^2+1)-x^{a-b}(x^3+x+1)=x(x^3-x^2+x-1).\end{multline}
The monomial $-x$ occurs in RHS, hence, one term of LHS has to be equal to $-x$. Thus, since $a>b$ and $c>b$ either 
\begin{eqnarray*}
\text{(case (3.5.1): } c-b=1, \text{ or} 
\text{(case (3.5.2): } a-b=1. 
\end{eqnarray*}
\item[(3.5.1)] If we replace $c$ by $b+1$ in (\ref{m:35}), we get
\begin{multline}\label{m:351}x^{\delta-b}(x^2+x+1)+x^{\epsilon-b}(x^4+x^2+1)-x^{a-b}(x^3+x+1)=x^2(x^2+2).\end{multline}
The monomial $2x^2$ occurs in RHS, therefore, two terms of LHS have to be equal to $x^2$. This follows that either
\begin{eqnarray*}
\text{(case (3.5.1.1): } \delta-b=1, \text{ or}
\text{(case (3.5.1.2): } \delta-b=2 
\end{eqnarray*} 
and either $\epsilon-b=0$, or $\epsilon-b=2$.

\item[(3.5.1.1)] By substituting $b+1$ for $\delta$ in (\ref{m:351}), we obtain $x^{\epsilon-b}(x^4+x^2+1)-x^{a-b}(x^3+x+1)=x(x^3-x^2+x-1).$
The monomial $-x$ occurs in RHS, thus, one term of LHS has to be equal to $-x$. Since $a \ne b$ we have $a-b=1$. This yields that $x^{\epsilon-b}(x^4+x^2+1)=x^2(2x^2-x+2).$ Obviously, LHS can be divided by $x^2+x+1$ and RHS cannot. Thus, this case cannot happen.
\item[(3.5.1.2)] If we replace $\delta$ by $b+2$ and $x$ by $t$ in (\ref{m:351}), we have $t^{\epsilon-b}(t^4+t^2+1)=t^2(1-t)$. Since $t^4+t^2+1=1-t$ we have $t^{\epsilon-b}=t^2$ and $\epsilon=b+2, a=b+3$. Thus, we get a pair of chromatically equivalent graphs $K_4(4,2,1,b+2,b+2,b) \cong K_4(4,2,1,b+1,b,b+3)$, where $b>1$.

\item[(3.5.2)] If we substitute $b+1$ for $a$ in (\ref{m:35}), we get 
$(x^{\delta-b}-x^{c-b})(x^2+x+1)+x^{\epsilon-b}(x^4+x^2+1)=x^2(2x^2-x+2).$ Evidently, LHS can be divided by $x^2+x+1$ and RHS cannot. In other words, this case cannot happen.

\item[(3.6)] By replacing $\epsilon$ by $b$ in (\ref{m:3}), we obtain
$(x^{\delta}-x^{c})(x^2+x+1)=(x^a-x^{\eta})(x^3+x+1).$ Obviously, since $\delta \ne c$ and $\eta \ne a$ RHS can be divided by $x^3+x+1$ and LHS cannot (roots of LHS are only either 1, or 0).

\item[(4)] Let us write down the equality.
\begin{multline*}
x^{\delta+5}+x^{\eta+5}+x^{\epsilon+4}+x^{\delta+2}+x^{\eta+2}+x^{\epsilon+3}-(x+1)(2x^{2}+x^{3}+x^{\delta}+x^{\epsilon}+x^{\eta})=\\ =x^{a+5}+x^{c+5}+x^{b+4}+x^{a+2}+x^{c+2}+x^{b+3}-(x+1)(2x^{2}+x^{3}+x^{a}+x^{b}+x^{c}).
\end{multline*}
Cancel the same items in both sides and divide by $x^2-1$.
\begin{multline}\label{m:4}
x^{\epsilon}(x^2+x+1)+x^{\eta}(x^2+1)(x+1)+x^{\delta}(x^2+1)(x+1)= \\ =x^{b}(x^2+x+1)+x^{a}(x^2+1)(x+1)+x^{c}(x^2+1)(x+1).
\end{multline}
Clearly, $\min\{\delta,\eta,\epsilon\}=\min\{a,b,c\}$. Therefore, w.l.o.g. either 
\begin{eqnarray*}
&&\text{case $(4.1)$: } \epsilon=b=1,  \\
&&\text{case $(4.2)$: } \epsilon=a=1,  \text{ or}\\
&&\text{case $(4.3)$: } \eta=c=1.\end{eqnarray*}

\item[(4.1)] By substituting 1 for $\epsilon$ and $b$ in (\ref{m:4}) and dividing by $x^2+1)(x+1)$, we have $x^{\eta}+x^{\delta}=x^{a}+x^{c}.$
Therefore, either $\eta=a, \delta=c$, or $\eta=c, \delta=a$. Both cases we obtain a pair of isomorphic graphs.

\item[(4.2)] If we substitute 1 for $\epsilon$ and $a$ in (\ref{m:4}) and divide by $x$, we get
$(x^{\eta-1}+x^{\delta-1}-x^{c-1})(x^3+x^2+x+1)-x^{b-1}(x^2+x+1)=x^3.$
W.l.o.g. $\eta \geq \delta$. Therefore, $\eta \geq c$ and $\eta \geq b-1$, otherwise monomial with maximal degree in LHS is negative. Thus, $\deg LHS=\eta-1+3$ and $\deg RHS=3$, hence, $\eta=1$. This is a contradiction (exactly one parameter of graph is equal to 1). The case $(4.2)$ cannot happen.

\item[(4.3)] If we replace $\eta$ and $c$ by 1 in (\ref{m:4}), we obtain
$(x^{\epsilon}-x^b)(x^2+x+1)=(x^a - x^{\delta})(x^2+1)(x+1)$. Whence, $(x^{\epsilon}-x^b)(x^3-1)=(x^a - x^{\delta})(x^4-1)$. Obviously, $\min\{\epsilon,b\}=\min\{a,\delta\}$. Since $\epsilon+\delta=a+b$ than $\epsilon=a$ and $\delta=b$. This case we obtain a pair of isomorphic graphs.

The proof of Theorem \ref{t:main} is complete.
\end{itemize}

\end{document}